\documentclass[preprint,12pt]{elsarticle}
\biboptions{sort&compress}

\usepackage{amssymb}
\usepackage{amsmath}
\usepackage{rotating}
\usepackage{pdflscape}
\usepackage{adjustbox}
\usepackage{makecell}
\usepackage{xcolor}
\usepackage{graphicx,caption,subcaption}
\usepackage{url}
\usepackage{hyperref}

\begin{document}

\begin{frontmatter}

\title{Patient Transport in Hospitals:\\ A Literature Review of Operations Research and Management Science Methods}

\author[inst1]{Tom Lorenz Klein}\corref{mycorrespondingauthor}
\ead{tom.klein@tum.de}
\cortext[mycorrespondingauthor]{Corresponding author}

\affiliation[inst1]{organization={Professorship of Optimization and Sustainable Decision Making,\\ Campus Straubing for Biotechnology and Sustainability,\\ Technical University of Munich},
            addressline={Am~Essigberg~3}, 
            city={Straubing},
            postcode={94315}, 
            country={Germany}}

\author[inst1,inst2]{Clemens Thielen}
\ead{clemens.thielen@tum.de}

\affiliation[inst2]{organization={Department of Mathematics, School of Computation, Information and Technology, Technical University of Munich},
            addressline={Boltzmannstr.~3}, 
            city={Garching bei München},
            postcode={85748}, 
            country={Germany}}

\begin{abstract}
 Most activities in hospitals require the presence of the patient. Delays in patient transport can disrupt operations, potentially resulting in idle staff, underutilized equipment, and postponed procedures, which in turn lead to lost revenue, unnecessary costs across many different areas and departments, and lower patient satisfaction. Consequently, patient transport planning is a central operational task in hospitals. This paper provides the first literature review of Operations Research and Management Science approaches for non-emergency, intra-hospital patient transport. We structure the different patient transport problems considered in the literature according to several main characteristics and introduce a five-field notation that allows for a concise representation of different problem variants. We then analyze the relevant literature with respect to different aspects related to the considered problem variant, the employed modeling and solution techniques, as well as the data used and the level of practical implementation achieved. Based on our literature analysis and semi-structured interviews with hospital practitioners, we compare current hospital practices and the existing literature, identify research gaps, and formulate an agenda for relevant future research.
\end{abstract}


\end{frontmatter}

\section{Introduction}\label{sec:intro}
Disrupted operations, delayed or postponed procedures, financial losses, and reduced quality of care are possible consequences of inefficient patient transport services in hospitals. While millions of patient transports are carried out in hospitals worldwide every year, poorly planned intra-hospital transports cause additional stress for patients and can lead to worse treatment outcomes and longer recovery times. They can also significantly increase the workload of transport staff, who travel up to 18 miles per day~\cite{salisburyjournal_2016}. In contrast, well-organized internal patient transports allow hospitals to better adhere to schedules for important operations such as surgeries or radiological examinations, resulting in less staff overtime and higher patient satisfaction. In addition, the loss of revenue from canceled surgeries is reduced, which in turn helps to ease the financial pressure that most hospitals face.

\medskip

Inefficient planning and coordination of predictable (i.e., non-emergency) intra-hospital patient transports (e.g., transfers between nursing wards or scheduled trips to the operating theater) not only lead to missed appointments, but also impact overall hospital operations and resource utilization both directly and indirectly. For example, late arrivals of patients in diagnostic departments such as radiology can disrupt the schedules of these departments, resulting in additional idle time for medical staff and undesired waiting time for patients. Moreover, also follow-up appointments might need to be rescheduled, which means that a single delay of a planned transport can lead to a series of downstream delays for surgeries and other procedures. This increases not only the scheduling complexity and uncertainty for all patient-related tasks, but also the required coordination between different hospital departments, which need to cooperate to fill gaps in their schedules and reorganize missed appointments. Ultimately, poor patient transport planning results in both lost hospital revenue due to suboptimal staff and equipment utilization, but also in lower patient satisfaction due to longer waiting times for appointments.

\medskip

While these considerations highlight the great importance of non-emer\-gen\-cy, intra-hospital patient transports, planning these transports is a major challenge due to the large number of location-related, time-related, and patient-related aspects that must be considered.

\medskip

Location-related aspects pertain to the origin, destination, and resulting route of a transport. Intra-hospital transports as considered here may involve routes inside of a hospital building, where patients are to be transported, e.g., from a nursing ward to the operating theater, or outdoor transports within a hospital complex or between hospitals. While transports within the same building usually mean that the patient is accompanied on foot or pushed in a bed or wheelchair, outdoor transports between different buildings often require using vehicles such as ambulances or minibuses. Depending on the specific situation, these vehicles might be allowed to transport several patients simultaneously, or be restricted to transporting only one patient at a time.

\medskip

For the non-emergency patient transports considered here, time-related aspects include time windows for pick-up and / or delivery of transported patients, which specify earliest and / or latest possible times for the start and / or the end of a transport. For example, an earliest possible pick-up time for a patient may result from the end of a previous appointment, and an earliest / latest delivery time may result from the scheduled procedure start at the destination department. While some of these times should be strictly adhered to whenever possible (e.g., arrival times of patients at the operating theater for surgery), others such as planned pick-up times of patients at their room on a nursing ward may offer more flexibility. Further important time-related aspects include working time regulations such as staff working hours and break times, or restrictions on the maximum duration of a transport that may result from a patient's medical condition.

\medskip

Patient-related aspects are diverse and range from specific transport needs to requiring additional medical assistance during transport. For transports within the same building, some patients can be accompanied on foot, while others depend on a bed or wheelchair, which must be provided in advance. In addition, certain patients may require continuous oxygen supply or special monitoring. Patients under the influence of medication may require medically trained transport staff, while infectious patients often lead to special isolation requirements, such as using protective clothing and disinfecting the transport equipment after each use. 
Overall, these factors lead to a diverse set of different transport scenarios, each of which entails specific logistical requirements. For example, a transport might involve accompanying a walking patient without any specific medical requirements from a nursing ward to the radiology department, while another transport might involve driving a patient requiring continuous oxygen supply and continuous monitoring by a medically trained person to another building within a hospital complex in an ambulance.

\medskip

The objectives pursued in patient transport planning are also diverse and lead to different planning approaches and strategies. 
One common objective is cost minimization, which focuses on reducing the transport costs resulting from personnel and equipment. In addition, there are strategies that focus on the efficient utilization of transport staff by minimizing both the walking / driving times and the waiting times of the transport personnel. A related objective focusing more on the patients is to reduce patient inconvenience caused by waiting and travel times. Other approaches take into account the ergonomic burden of transport personnel in order to improve their working conditions. When considering multiple objectives simultaneously, conflicts may arise and compromise solutions might have to be found. For example, reducing patient waiting times by utilizing additional transport personnel can increase the total transport costs, while choosing a transport schedule that leads to a high staff and equipment utilization may lead to longer waiting times for patients.

\medskip

In addition to the aspects mentioned above and the diverse, complex requirements and conflicting objectives, intra-hospital patient transport planning is often complicated by uncertainty - even when non-emergency transports are considered. For instance, additional transport requests may arise during the day due to newly-admitted patients or additional examinations or treatments that are requested for patients that have already been admitted before. Moreover, transports might be postponed due to delays in other patient-related tasks, and the time required to perform a transport may change due to unforeseen factors such as a blocked elevator. This makes the practical problem a dynamic problem that involves various sources of uncertainty.

\medskip

All of these factors make non-emergency, intra-hospital patient transport planning an important and challenging task, for which a variety of Operations Research and Management Science (OR/MS) approaches have been developed.
This paper provides the first literature review of these approaches. Moreover, we identify and highlight differences in the aspects of patient transport modeled in the relevant literature as well as existing research gaps. To align this analysis with requirements from hospital practice and derive a practically relevant research agenda, we conducted semi-structured interviews with practitioners from five hospitals in Germany and the Netherlands and compare the resulting insights with the current state of the literature.

\medskip

The remainder of this paper is structured as follows: Section~\ref{sec:search} presents the methodological approach to the literature search, which led to the identification of the relevant scientific contributions. Afterwards, Section~\ref{sec:notation} outlines the existing modeling approaches found in the identified literature. Moreover, a comprehensive classification of different variants of patient transport problems by means of a five-field notation is developed. We use this notation in Section~\ref{sec:statistics} as a basis for classifying the relevant publications from the literature search. The section then provides a detailed discussion of key characteristics of patient transport problems considered in the examined literature and analyzes further relevant aspects such as modeling and solution methods. Moreover, cross comparisons are made between the considered characteristics. Section~\ref{sec:gap} then identifies research gaps and formulates an agenda for future research based on our literature analysis and interviews with hospital practitioners. Finally, Section~\ref{sec:conclusion} concludes the paper with a summary and discussion of the obtained insights.

\section{Literature search methodology} \label{sec:search}
\begin{figure}[h]
    \centering
    \includegraphics[width=0.70\textwidth]{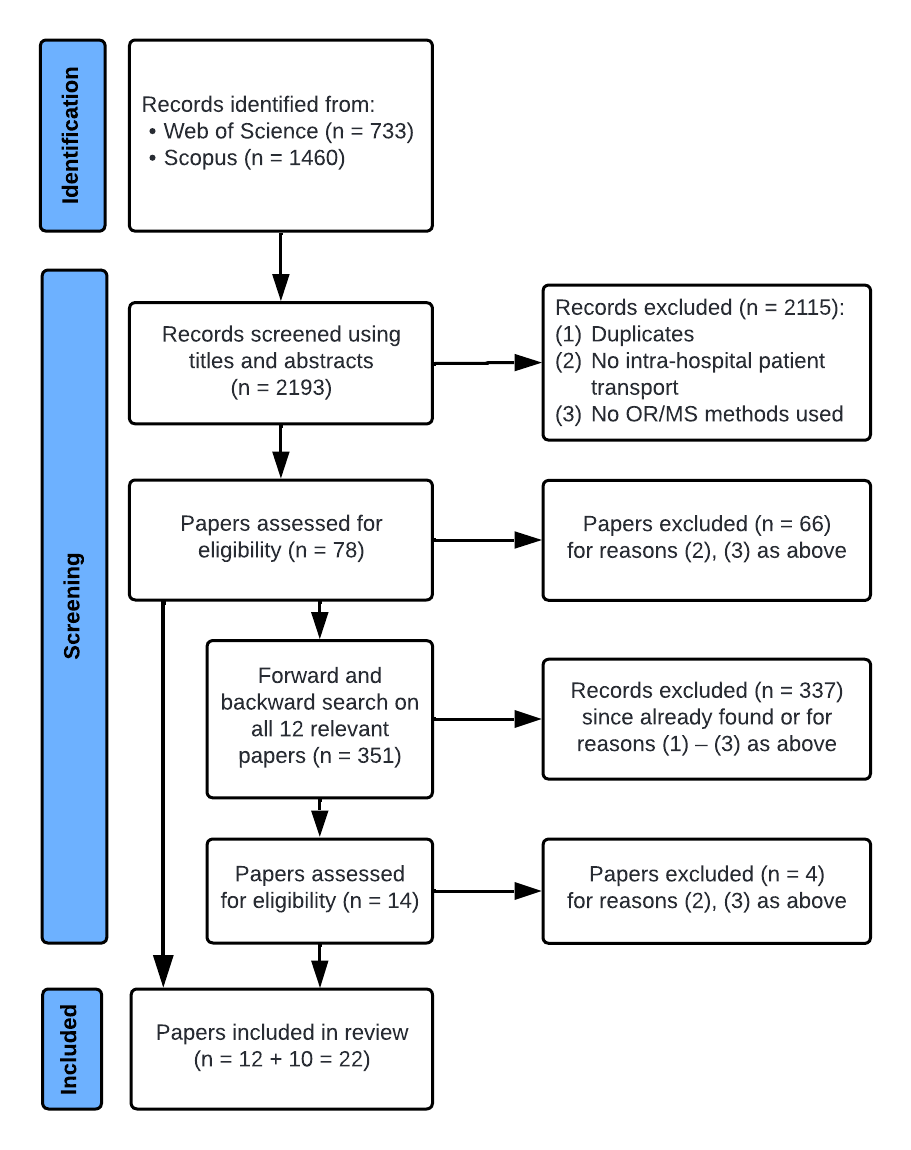}
    \caption{PRISMA 2020 flow diagram illustrating our literature search methodology.}
    \label{fig:prisma_chart}
\end{figure}

In order to provide a sound and comprehensive basis for our review, an in-depth literature search was carried out in October 2024 using the databases Scopus and Web of Science within journals classified as ``Operations Research \& Management Science'' according to their All Science Journal Classification Code (ASCJ) in Scopus or their Web of Science Category or research area (or both) in Web of Science. Our literature search methodology follows the PRISMA guidelines (\url{www.prisma-statement.org}) and is illustrated using a PRISMA 2020 flow diagram~\cite{Page_prisma} in Figure~\ref{fig:prisma_chart}.

To ensure that no relevant literature is overlooked, the number of required search terms was deliberately restricted to a minimum. The requirement was that at least one term from each column in Table~\ref{tab:keywords} appears in the title, the abstract, or the author keywords. Here, the first column ``Hospital terms'' ensures the connection of the publications to hospitals, while the second column ``Transport terms'' ensures a transport focus. Note that, since patient transport problems are often modeled as vehicle routing problems (VRPs) and, in particular, as Dial-a-Ride Problems (DARPs), corresponding search terms have been included in the transport terms column. 

\medskip
\begin{table}[h]
\centering

    \begin{tabular}{c|c}
        \textbf{Hospital terms}  & \textbf{Transport terms} \\
        \hline
        hospital\$ 	&	dial*    \\
        clinic* 	&	ride*\\
        infirmary 	&	vehicle*\\
        infirmaries 	&	rout*\\
        “medical school\$” 	&	transp*\\
        “college\$ of medicine” 	&	DARP*\\
        “department\$” 	&	VRP*\\
        surger* 	&	\\
        surgeo*	&	\\
        surgical	&	\\
        patient\$ &	\\

    \end{tabular}
    \centering
    \caption{Hospital and transport search terms. The wildcard ``\$'' stands for at most one character, while the wildcard ``$\ast$'' stands for any group of characters (including no characters). The search terms in each column are combined with an \emph{OR} operator, while an \emph{AND} operator is used to link hospital and transport terms (i.e., at least one hospital term and at least one transport term must be found). }
    \label{tab:keywords}
\end{table}

\medskip

The search yielded 1460 results in the Scopus database and 733 in the Web of Science database. Only research items (journal and conference papers, PhD theses, books, and book chapters) in English were reviewed. Titles and abstracts were carefully examined to identify studies specifically addressing non-emergency patient transport within hospitals,
such as planned transfers from nursing wards to radiology or to the operating theater. However, considering that some hospitals consist of multiple buildings or even different locations within a geographical region, papers that focus on patient transport between hospital buildings (of the same hospital or of different hospitals) were also included.

\medskip
After the title and abstract screening described above, a total of 78 potentially relevant papers were left. The full texts of these papers were then read in detail. Here, only papers that focus on intra-hospital patient transport and use methods from Operations Research and Management Science were selected as relevant. 
Afterwards, a forward and backward search was performed on the selected papers to avoid overlooking relevant papers, e.g., papers not published in journals classified as OR/MS.

\medskip

Overall, using the process described above, we identified 22~papers applying OR/MS approaches to patient transport within hospitals or between hospital complexes, which were included in the review. Table~\ref{tab:paper_overview} provides an overview of these papers, while Table~\ref{tab:paper_five_field} in Section~\ref{sec:statistics} classifies each paper according to the five-field notation that will be introduced in the following section.

\begin{table*}
    \centering
    \rotatebox{90}{
    \begin{minipage}{0.94\textheight}
    \caption[Identified relevant papers]{Identified relevant papers sorted by year of publication, including publication outlet and investigated region. Entries marked with $\ast$ represent cases where the region could not be directly identified and instead reflects the region of the majority of the authors.
    }
    \label{tab:paper_overview}
\scalebox{0.75}{
    \begin{tabular}{|l|l|l|l|}
    \hline

         Publication & Publication outlet & Year & Region \\
      
    \hline

         Nickel et al.~\cite{nickel_2000a} & Operations Research Proceedings 1999 & 2000 & *Germany \\
         Kallrath~\cite{kallrath_2005a} & Book chapter & 2005 & Germany \\
         Melachrinoudis et al.~\cite{melachrinoudis_2005a} & Computers \& Operations Research & 2007 & USA \\
         Fiegl et al.~\cite{fiegl_2009a} & Journal of Biomedical Informatics & 2009 & Austria \\
         Hanne et al.~\cite{hanne_2009a} & Interfaces & 2009 & Germany \\
         Beaudry et al.~\cite{beaudry_2010a} & OR Spectrum & 2010 & Germany \\
         Kergosien et al.~\cite{kergosien_2011a} & European Journal of Operational Research & 2011 & France \\
         Turan et al.~\cite{turan_2011a} & IEEE International Symposium on Logistics and Industrial Informatics & 2011 & *Austria \\
         Kergosien et al.~\cite{kergosien_2014a} & International Conference on Health Care Systems Engineering & 2014 & Canada \\
         Schmid et al.~\cite{schmid_2014a} & Transportation Science & 2014 & Austria \\
         von Elmbach et al.~\cite{elmbach_2015a} & IIE Transactions on Healthcare Systems Engineering & 2015 & Germany \\
         Zhang et al.~\cite{zhang_2015a} & Omega & 2015 & Hong Kong \\
         Detti et al.~\cite{detti_2017a} & Omega & 2017 & Italy \\
         Vancroonenburg et al.~\cite{vancroonenburg_2018a} & International Conference on the Practice and Theory of Automated Timetabling & 2016 & Belgium \\
         Séguin et al.~\cite{seguin_2019a} & Operations Research for Health Care & 2019 & Canada \\
         von Elmbach et al.~\cite{elmbach_2019a} & European Journal of Operational Research & 2019 & *Germany \\
         van den Berg et al.~\cite{berg_2019a} & Transportation Science & 2019 & Netherlands \\
         Xiao et al.~\cite{xiao_2022a} & Computers \& Operations Research & 2022 & China \\
         Nasira et al.~\cite{nasir_2022a} & Journal of Transport \& Health & 2022 & Hong Kong \\
         Maka et al.~\cite{maka_2022a} & LogForum & 2022 & Thailand \\
         Kergosien et al.~\cite{kergosien_2023a} & International Journal of Production Research & 2023 & Canada \\
         Bärmann et al.~\cite{Baermann2024} & OR Spectrum & 2024 & Germany \\
     
    \hline
    \end{tabular}
   }
    \end{minipage}
    }
\end{table*}

\section{Modeling approaches and classification}\label{sec:notation}
In this section, we further elaborate on the types of patient transport problems studied in the literature. We begin by detailing different modeling approaches in Subsection~\ref{sec:darp}, placing patient transport within the context of related OR/MS problems. Following this, Subsection~\ref{sec:classification} introduces a five-field notation for classifying patient transport problem variants, providing a structured overview of their distinguishing characteristics as identified in the literature.

\medskip

\subsection{Modeling approaches and relations to vehicle routing problems}\label{sec:darp}
Existing research highlights several modeling frameworks for patient transport, with most studies characterizing the problem as a variant of the Vehicle Routing Problem (VRP) or a closely related problem~\cite{toth_vigo_2002a}. Notably, Séguin et al.~\cite{seguin_2019a} and Maka et al.~\cite{maka_2022a} model the patient transport problem as an assignment and resource allocation problem, respectively, while Fiegl et al.~\cite{fiegl_2009a} and von~Elmbach et al.~\cite{elmbach_2015a} represent it as a scheduling problem. However, most studies employ VRP-based models, specifically focusing on the Capacitated VRP (CVRP) for handling vehicle capacities, the VRP with Pick-up and Delivery (PDVRP) to model location-based patient transfers, and the VRP with Time Windows (VRPTW) for managing scheduling constraints, such as patient pick-up and delivery time windows.

\medskip

The Dial-a-Ride Problem (DARP) (see, e.g.,~\cite{cordeau_2007a}) is the most frequently applied model in the literature since it allows modelling a wide range of aspects such as patient travel and waiting times, isolated transports, the adherence to working hours of transport staff, skill or equipment requirements of transports, and patient inconvenience. Here, patient inconvenience refers to any discomfort or additional burden experienced by patients due to delays, prolonged transport times, or extended waiting periods before transport. For example, elderly or severely ill patients may experience fatigue or pain if the travel time is extended beyond a certain limit. Similarly, long waiting times before a scheduled transport can cause discomfort for patients and increase anxiety or stress. This holds particularly in cases where the patient is already in a weakened state or anxious about an upcoming procedure, or when the patient waits for a transport back to their room after an exhausting appointment.

\medskip

\subsection{Classification of patient transport problems}\label{sec:classification}
Although many studies model patient transport problems as variants of the DARP, the specific features considered vary significantly. To facilitate a more granular understanding of these variants, we introduce a five-field notation~$\alpha \,\vert\, \beta \,\vert\, \gamma \,\vert\, \delta \,\vert\,\epsilon$ for classifying patient transport problems, drawing inspiration from the well-known notation for scheduling problems by Graham et al.~\cite{graham_1979a}. Our notation includes five fields: fleet characteristics~($\alpha$), depot characteristics~($\beta$), constraints~($\gamma$), objectives~($\delta$), and uncertainty~($\epsilon$), which are further detailed in Table~\ref{tab:five_field_overview} and the subsequent sections.

\begin{table}

    \centering
    \begin{tabular}{|c|l|c|l|}
    \hline
        \multicolumn{2}{|c|}{Fleet characteristics~$\alpha$}  & \multicolumn{2}{c|}{Depot characteristics~$\beta$} \\
        \hline
         $F_{\text{MULT}}$& multiple vehicle types  & $D_{\text{MULT}}$ & multiple depots \\
         $F_{N}$& $N$ vehicle types  & $D_{N}$ & $N$ depots \\
          $F_1$ & single vehicle type &  $D_1$& single depot \\
          $F_\emptyset$ & vehicle-free  & $D_\emptyset$ &no depots  \\
         \hline
    \end{tabular}\\
   \bigskip
    \centering
    \begin{adjustbox}{width=\textwidth,keepaspectratio=true}
    \begin{tabular}{|c|l|c|l|c|l|c|l|} 
    \hline
    \multicolumn{8}{|c|}{Constraints~$\gamma$} \\
    \hline
    \multicolumn{2}{|c|}{Capacity} & \multicolumn{2}{c|}{Equipment} & 
    \multicolumn{2}{c|}{Priority} &
    \multicolumn{2}{c|}{Skill}  \\
    \hline
    $CAP_1$&unit & $EQ$ & additional & $\text{PRIO}$ & priority & $SK$& skill  \\
    $CAP_{\text{UNIF}}$& uniform & & equipment & & levels & & requirements    \\
    $CAP_{\text{MULT}}$& multi& & & & & &  \\
    \hline
    \end{tabular}
    \end{adjustbox}
    \\ 
    \bigskip
    \begin{tabular}{|c|l|c|l|c|l|} 
    \hline
    \multicolumn{6}{|c|}{Constraints~$\gamma$} \\
    \hline
    \multicolumn{2}{|c|}{Isolation} & \multicolumn{2}{c|}{Time windows} &  \multicolumn{2}{c|}{Time} \\
    \hline
    $ISO$& isolated & $TW_{HARD}$ & hard   &$T_{S}$ & staff  \\
    & transports& $TW_{SOFT}$& soft & $T_{L}$ & loading  \\
    & & $TW_{MIX}$ & mixed  &$T_{TT}$ & travel time \\
    \hline
    \end{tabular}
    \\
    \bigskip
    \begin{adjustbox}{width=\textwidth,keepaspectratio=true}
    \begin{tabular}{|c|l|c|l|c|l|c|l|c|l|} 
    \hline
    \multicolumn{10}{|c|}{Problem objective(s)~$\delta$} \\
    \hline
    \multicolumn{2}{|c|}{Cost} & \multicolumn{2}{c|}{Waiting time} & \multicolumn{2}{c|}{Travel time} & \multicolumn{2}{c|}{Ergonomic burden} & \multicolumn{2}{c|}{Flow time}  \\
    \hline
        $C$& cost &$WT_{V}$ & vehicles &$TT_{V}$ & vehicles &$\max B$ & maximum&$\sum w_jF_j$ & total  \\
        & &$WT_{P}$ & patients  &$TT_{P}$ & patients & &ergonomic & & weighted\\
        & & & &$TT_{E}$ & empty & & burden& & flow time \\
    \hline
    \end{tabular}
    \end{adjustbox}
    \\
    \bigskip
    \begin{tabular}{|c|l|}
    \hline
        \multicolumn{2}{|c|}{Uncertainty~$\epsilon$}  \\
        \hline
         $R$& uncertain requests  \\
         $TT$& uncertain travel times \\
         \hline
    \end{tabular}
    \caption{Overview of the five-field notation $\alpha \,\vert\, \beta \,\vert\, \gamma \,\vert\, \delta \,\vert\,\epsilon$ with a description of the possible entries for each field.}
    \label{tab:five_field_overview}
\end{table}

\subsubsection{Fleet characteristics ($\alpha$)}\label{sec:fleetchar}
The first field~$\alpha$ of our five-field notation classifies the fleet characteristics. Here, the term \emph{fleet} refers to the group of people or manned vehicles responsible for patient transports. In many of the examined papers, patient transports inside hospitals are exclusively carried out by transport personnel accompanying a walking patient or pushing a patient in a wheelchair or bed, which we refer to as a \emph{vehicle-free fleet} denoted by $\alpha \equiv F_\emptyset$. In cases where vehicles are used, the characteristics in which the fleet can differ are vehicle capacity, speed, and range, the ability to transport equipment, and the option to execute isolated transports. A \emph{multiple-vehicle-types fleet} or \emph{heterogeneous fleet} that contains at least two vehicle types that differ in at least one of these characteristics is denoted by $\alpha \equiv F_{\text{MULT}}$. To specify a fleet with a specific number~$N$ of vehicle types, the notation~$\alpha \equiv F_{N}$ is used. In particular, a \emph{single-vehicle-type fleet} or \emph{homogeneous fleet} of identical vehicles that do not differ in any of the mentioned characteristics is denoted by~$\alpha \equiv F_1$. 

\subsubsection{Depot characteristics ($\beta$)}\label{sec:depotchar}
The second field~$\beta$ characterizes the modeling of depots (hubs, assembly points) for the transport staff or vehicles. In case \emph{multiple depots} are considered, we use the notation~$\beta \equiv D_{\text{MULT}}$. Moreover, we use~$\beta \equiv D_{N}$ for problems with a specific number~$N$ of depots. In particular, some relevant papers use a \emph{single-depot} approach, which we denote by~$\beta \equiv D_1$. Other problems do not use depots at all, which we denote by $\beta \equiv D_\emptyset$. These problems assume that the transport staff (e.g., nurses) work on the nursing wards that are also the pick-up locations of the patients they transport, which means that no travel is required in order to reach the pick-up location of a transport. 

\subsubsection{Constraints ($\gamma$)}\label{sec:conchar}
The third field~$\gamma$ indicates the presence of several kinds of important constraints that can be considered. These include capacity constraints, equipment constraints, constraints concerning transport priorities or skill requirements for transport personnel, as well as constraints regarding isolated transports, time windows, and other time-related constraints. If more than one of these is considered, the third field will have multiple entries, which are sorted in the same order as in Table~\ref{tab:five_field_overview} and separated by semicolons.

\medskip 
The first entry in the third field specifies the capacity constraints. Here, the entry $CAP_1$ indicates \emph{unit capacity}, which means that each transporter or transport vehicle can only transport one unit (whether bed, wheelchair, or walking patient) at a time. Unit capacity modeling is often used for the transport of disabled patients and for ambulance transport (with some exceptions, where ambulances can transport multiple patients simultaneously~\cite{beaudry_2010a,hanne_2009a}). \emph{Uniform capacity}, indicated by $CAP_{\text{UNIF}}$, describes that all transporters or transport vehicles have a uniform capacity strictly larger than one. Finally, we indicate \emph{multi-capacity problems} by  $CAP_{\text{MULT}}$. In this case, the fleet includes at least two vehicles with different capacities.

\medskip

In addition, challenges arise during patient transport when special equipment such as oxygen supply or monitoring equipment is required. Due to this specialized equipment or other factors such as the patient being under the influence of medication, different skill levels of transport personnel may also be required. Some problems also include priority levels that affect the sequence in which the requests are processed. Furthermore, in certain cases, transport requests may need to be handled separately depending on the patient's condition or other requirements. To represent the consideration of transport requests requiring \emph{special equipment}, we use the symbol~$EQ$ in the third field. \emph{Skill requirements}, which need to involve at least two levels, are denoted by~$SK$. For problems that involve at least two priority levels of transports or isolated transports, we use the symbols~$PRIO$ and~$ISO$, respectively.

\medskip

Moreover, patient transport problems include different approaches to \emph{time window} constraints, which specify earliest and / or latest possible times for the start (pick-up) and / or the end (delivery) of each transport. \emph{Hard time windows} are denoted by~$TW_{\text{HARD}}$, while \emph{soft time windows} for which penalties are imposed for early and / or late pick-ups and / or deliveries are denoted by~$TW_{\text{SOFT}}$. In case of problems that require strict adherence to some bounds of the imposed time windows but allow others to be violated (e.g., latest possible end times of transports cannot be violated, but earliest possible end times can), we refer to \emph{mixed time windows}, which are denoted by~$TW_{\text{MIXED}}$.

\medskip

Besides time windows, we also indicate the presence of other kinds of \emph{time-related constraints} in the third field. These include \emph{staff time constraints} on the compliance with working times and the consideration of break times of the transport staff, which we denote by~$T_{S}$. Likewise, some problem definitions consider times for disinfecting transport vehicles or other times for preparation and post-processing, as well as the times for loading and unloading of patients. We jointly refer to these as \emph{loading time constraints} and denote them by~$T_{L}$. In cases where multiple patients are transported together on a single tour or when transporting severely ill patients, a maximum \emph{travel time} ($T_{TT}$) can be defined for each individual patient, which limits the time between this particular patient's pick-up and delivery (i.e., it applies to the individual patient's travel time, and not to the entire tour). These time-related constraints are not exclusive, i.e., staff time constraints, loading time constraints, and travel time constraints may all be considered together, which is then indicated by multiple indices at the letter~$T$. For example, a problem that considers all of these constraints will have an entry $T_{S,L,TT}$ in the third-field.

\subsubsection{Problem objective(s) ($\delta$)}\label{sec:objectchar}

 The fourth field describes the \emph{problem objective(s)}. Here, the first important distinction is between single-objective problems, multi-objective problems, and problems that use a weighted sum objective. Single-objective problems use only one of the minimization objectives described below. In multi-objective problems, where several objectives are considered separately, the different minimization objectives are separated using semicolons. To show that a problem uses a weighted sum objective, where a single minimization objective is obtained as a nonnegative weighted sum of several objectives, the minimization objectives are connected by a plus (``$+$''). Note that, in cases where the different objectives have different units of measurement (e.g., when forming a weighted sum of a cost and a time objective), the weights must, in particular, convert all objectives in the weighted sum to a common unit.

\medskip

Another frequently considered objective is cost minimization. There are many different and complex combinations of costs found in the literature, so we simply indicate \emph{cost minimization} by an entry~$C$ in the fourth field. Examples of cost objectives found in the literature include staff costs, transport costs (consisting of fuel costs, maintenance costs, cleaning costs, and general vehicle costs), and many more. 

\medskip

A further objective is the minimization of \emph{waiting times}. Here, we distinguish between the minimization of waiting times for \emph{vehicles}, denoted by~$WT_{V}$, and the minimization of waiting times for \emph{patients}, denoted by~$WT_{P}$. Waiting times for vehicles arise if vehicles are already on site earlier than a specified pick-up time or if it is required that a patient may not be unloaded before a specified unloading time. Waiting times for patients arise from the difference between the earliest possible pick-up time and the actual pick-up time, as well as in the event of earliness at the delivery location. 

\medskip

A related type of objective considered is the minimization of \emph{travel times}. Here, we distinguish the minimization of \emph{travel times of vehicles} ($TT_{V}$), \emph{travel times of patients} ($TT_{P}$), and \emph{empty travel times of vehicles} ($TT_{E}$). Minimizing travel times of vehicles ($TT_{V}$) might result in longer travel times for patients, as this approach focuses on reducing the duration of vehicle usage and optimizing operational efficiency. This is particularly important in scenarios with limited vehicle availability or high operational costs of vehicles. Conversely, minimizing patient travel time ($TT_{P}$) often involves choosing direct routes instead of longer tours with intermediate stops, prioritizing patient comfort and well-being. 
The empty travel time of vehicles ($TT_{E}$) is also a key factor. It refers to the time when vehicles are traveling without patients, either to reach a patient or to return (e.g., to the depot) after completing a transport. Minimizing empty travel times of vehicles ($TT_{E}$) is aimed at reducing operational costs and improving overall service efficiency. Thus, striking a balance between minimizing $TT_{V}$, $TT_{P}$, and $TT_{E}$ is essential. This involves a trade-off between operational efficiency and patient needs, ensuring that the transport system is both cost-effective and patient-centered. 

\medskip

Two further objectives considered are the minimization of the \emph{maximum ergonomic burden} of the transport staff, denoted by~$\max B$, and the minimization of the \emph{total weighted flow time}, denoted by $\sum w_jF_j$. Here, the ergonomic burden depends on several factors such the tour length, the transported patient's weight, or how the patient is transported (e.g, in a bed or in a wheelchair). For the total weighted flow time -- as is common for jobs in scheduling problems -- the flow time~$F_j$ of a transport request~$j$ is its completion time minus its release date (earliest possible start time).

\subsubsection{Uncertainty ($\epsilon$)}\label{sec:uncchar}

The fifth field~$\epsilon$ addresses the concept of \emph{uncertainty} in patient transport problems. Here, we distinguish between two main forms of uncertainty. First, there may be \emph{uncertain requests}, denoted by an entry~$R$. This refers to situations where not all transport requests are known in advance, which can happen, e.g., when additional transport requests are triggered due to additional examinations or procedures being requested for a patient during a doctor's visit. Second, \emph{uncertain travel times} are denoted by~$TT$, which refers to the case where the required travel time of a transport request is not precisely known before executing it (e.g., due to varying waiting times for elevators). Note that both forms of uncertainty may be considered either individually or in combination. For instance, a problem may include uncertain travel times even if all requests are known in advance, and a problem with uncertain requests may have either known travel times (if the exact travel time becomes known at the time a transport request is received) or uncertain travel times (if the travel time remains uncertain until the actual execution of the transport request).
\section{Literature classification and analysis}\label{sec:statistics}
In this section, we first demonstrate the versatility of the five-field notation introduced in Section~\ref{sec:notation} by classifying the problems described in the relevant articles identified in our literature search. We then present a detailed analysis of the problems considered in these articles according to various criteria.

\begin{table*}
    \rotatebox{90}{
    \begin{minipage}{0.92\textheight}
    \caption[Classification of patient transport problems.]{Classification of patient transport problems in the papers identified in our literature search (in chronological order) using the five-field notation introduced in Table~\ref{tab:five_field_overview}.}
    \label{tab:paper_five_field}
    \scalebox{0.75}{
    \begin{tabular}{|l|c|c|l|l|c|}
    \hline
         Publication & $\alpha$ & $\beta$ & $\gamma$ & $\delta$ & $\epsilon$ \\
    \hline
         Nickel et al.~\cite{nickel_2000a} & $F_\emptyset$& $D_1$ & $CAP_1;PRIO$ & $C;TT_{E}$& $R$  \\
         Kallrath~\cite{kallrath_2005a}& $F_{\text{MULT}}$ & $D_{\text{MULT}}$ & $CAP_{\text{MULT}};EQ;ISO;TW_{\text{SOFT}};T_{TT}$ & $C+WT_{V,P}+ TT_{V,P}$ & $R$ \\
         Melachrinoudis et al.~\cite{melachrinoudis_2005a} & $F_{\text{MULT}}$ & $D_{\text{MULT}}$& $CAP_{\text{MULT}};TW_{\text{SOFT}}$ & $C+WT_{P}+TT_{P}$ & \\
          Fiegl et al.~\cite{fiegl_2009a}&$F_\emptyset$&$D_\emptyset$&$CAP_{\text{MULT}};PRIO;TW_{\text{SOFT}}$& $\sum w_jF_j$  & $R$\\
         Hanne et al.~\cite{hanne_2009a}& $F_{\text{MULT}}$& $D_{\text{MULT}}$ & $CAP_{\text{MULT}};EQ;ISO;T_{S,L,TT}$ & $C+WT_{V,P}+TT_{V,P}$ & $R,TT$  \\
         Beaudry et al.~\cite{beaudry_2010a}& $F_{\text{MULT}}$& $D_{\text{MULT}}$ & $CAP_{\text{MULT}};EQ;PRIO;SK;ISO;TW_{\text{MIX}};T_{S,L,TT}$  & $C+WT_{V,P}+TT_{V,P,E}$ & $R,TT$  \\
         Kergosien et al.~\cite{kergosien_2011a} & $F_{\text{MULT}}$& $D_{\text{MULT}}$ & $CAP_1;EQ;PRIO;ISO;TW_{\text{HARD}};T_{S}$ & $C$ & $R$ \\
         Turan et al.~\cite{turan_2011a}& $F_\emptyset$& $D_{\text{MULT}}$ & $CAP_1;PRIO;TW_{\text{MIX}};T_{L}$ & $C+WT_{V,P}+TT_{V,E}$ & \\
         Kergosien et al.~\cite{kergosien_2014a} & $F_1$ & $D_{\text{MULT}}$ & $CAP_1;PRIO;TW_{\text{MIX}};T_{S}$ & $C; WT_{P};TT_{E}$ &$R$ \\
         Schmid et al.~\cite{schmid_2014a} & $F_\emptyset$& $D_{\text{MULT}}$ & $CAP_1;TW_{\text{MIX}}$ & $C+WT_{V,P}+TT_{E}$ & \\
         von Elmbach et al.~\cite{elmbach_2015a}& $F_\emptyset$ & $D_{\text{MULT}}$ & $CAP_1;PRIO;TW_{\text{MIX}}$ & $\max B$ &  \\
         Zhang et al.~\cite{zhang_2015a}& $F_1$ & $D_{\text{MULT}}$ & $CAP_{\text{MULT}};SK;TW_{\text{HARD}}; T_{S,L,TT}$ & $C+WT_{V,P}+TT_{V,P}$ & \\
         Detti et al.~\cite{detti_2017a}& $F_{\text{MULT}}$& $D_{\text{MULT}}$ & $CAP_{\text{MULT}};ISO;TW_{\text{MIX}};T_{L,TT}$ & $C+WT_{V,P}+TT_{V,P,E}$ &  \\
         Vancroonenburg et al.~\cite{vancroonenburg_2018a} & $F_\emptyset$ & $D_{\text{MULT}}$ & $CAP_1;SK; TW_{\text{MIX}}; T_{S}$ & $TT_{V}$ &$R$ \\
         Séguin et al.~\cite{seguin_2019a}& $F_\emptyset$ & $D_\emptyset$ & $T_{S}$ & $C$ & \\
         von Elmbach et al.~\cite{elmbach_2019a} & $F_\emptyset$ & $D_{\text{MULT}}$ & $CAP_1;TW_{\text{SOFT}};T_{L}$ & $\max B$&  \\
         Van den Berg et al.~\cite{berg_2019a} & $F_{\text{MULT}}$& $D_{\text{MULT}}$ & $CAP_1;PRIO;SK;TW_{\text{MIX}};T_{S}$ & $TT_{V}$ & $R,TT$ \\
         Xiao et al.~\cite{xiao_2022a}& $F_\emptyset$ & $D_\emptyset$ & $CAP_1;TW_{\text{SOFT}}$ & $WT_{P}+TT_{P}$ &  \\
         Nasira et al.~\cite{nasir_2022a} & $F_1$ & $D_1$ & $CAP_{\text{MULT}};TW_{\text{SOFT}}; T_{L}$ & $C+WT_{V,P}+TT_{V,P,E}$ & \\
         Maka et al.~\cite{maka_2022a}& $F_\emptyset$ & $D_{\text{MULT}}$ & $CAP_1;TW_{\text{HARD}};T_{S,L,TT}$ & $C+WT_{P}+TT_{V,P}$ &  \\
         Kergosien et al.~\cite{kergosien_2023a}& $F_1$ & $D_1$ & $CAP_1;SK;TW_{\text{SOFT}};T_{S,L}$ & $C+TT_{V}$ & $R$  \\
         Bärmann et al.~\cite{Baermann2024} & $F_\emptyset$ & $D_\emptyset$ & $CAP_1;PRIO;SK;ISO;TW_{\text{MIX}};T_{S}$ & $TT_{E}$ &$R$  \\
    \hline
    \end{tabular}
    }
    \end{minipage}
    }
\end{table*}

\subsection{Transport origins, destinations, and items}

We begin by examining the transport origins and destinations, as well as the transported items. As already stated in Section~\ref{sec:search}, the focus of this review is on intra-hospital patient transport, i.e., transport of patients within a hospital building or complex or between different hospitals. However, some of the selected papers also consider other kinds of transports. Therefore, Figures~\ref{fig:locations} and~\ref{fig:goods} show the distribution of the transport origins / destinations and the transported items, respectively, in the relevant literature. 

\begin{figure}[t]
  \begin{subfigure}[t]{\textwidth}
    \centering
    \includegraphics[width=\linewidth]{locations.png}
    \caption{Distribution of origin-destination pairs.}
    \label{fig:locations}
    \footnotesize
    \vspace{0.5em}
    \begin{tabular}{ll}
        Inside hospitals & \cite{nickel_2000a,kallrath_2005a,melachrinoudis_2005a,fiegl_2009a,hanne_2009a,beaudry_2010a,turan_2011a,schmid_2014a,elmbach_2015a,vancroonenburg_2018a,seguin_2019a,elmbach_2019a,xiao_2022a,maka_2022a,Baermann2024}  \\
        Hospital to hospital & \cite{melachrinoudis_2005a,hanne_2009a,kergosien_2011a,kergosien_2014a,zhang_2015a,detti_2017a,berg_2019a,nasir_2022a,kergosien_2023a} \\
        Hospital to home & \cite{melachrinoudis_2005a,kergosien_2014a,detti_2017a,berg_2019a,nasir_2022a,kergosien_2023a} \\
        Accident location to hospital & \cite{kergosien_2014a,berg_2019a,kergosien_2023a} \\
    \end{tabular}
    \label{tab:origins}
  \end{subfigure}
  
  \hfill
  
  \begin{subfigure}[t]{\textwidth}
    \centering
    \includegraphics[width=\linewidth]{goods.png}
    \caption{Distribution of transported items.}
    \label{fig:goods}
    \footnotesize
    \vspace{0.5em}
    \begin{tabular}{ll}
        Patient & \cite{nickel_2000a,kallrath_2005a,melachrinoudis_2005a,fiegl_2009a,hanne_2009a,beaudry_2010a,kergosien_2011a,turan_2011a,kergosien_2014a,schmid_2014a,elmbach_2015a,zhang_2015a,detti_2017a,vancroonenburg_2018a,seguin_2019a,elmbach_2019a,berg_2019a,xiao_2022a,nasir_2022a,maka_2022a,kergosien_2023a,Baermann2024}  \\
        Material & \cite{kallrath_2005a,melachrinoudis_2005a,vancroonenburg_2018a,Baermann2024} \\
        Emergency patient & \cite{kergosien_2014a,berg_2019a,kergosien_2023a} \\
    \end{tabular}
    
  \end{subfigure}
  
  \caption{Distribution of different transport origins/destinations and transported items in the examined papers.}
  \label{fig:loc-goods}
\end{figure}

\medskip

\enlargethispage{\baselineskip}

In Figure~\ref{fig:locations}, transports inside hospital buildings or between buildings of a hospital complex are denoted by \emph{Inside hospitals}, while transports of patients between hospitals (e.g., planned transfers to a different hospital) are referred to as \emph{Hospital to hospital}. Transports of patients from home to hospital or discharges from hospital to home are denoted by \emph{Hospital to home} and transports from an accident site to the hospital, which are usually emergency transports, are denoted by \emph{Accident location to hospital}.
As can be seen from Figure~\ref{fig:locations}, each paper included in our review considers either transports inside hospitals or hospital to hospital transports (or both), which reflects the requirements of our literature selection. However, transports between a patient's home and a hospital are also considered in a relevant number of papers, which indicates that these transport problems are closely linked to transports inside and between hospitals. Transports from an accident location to a hospital are examined particularly in~\cite{kergosien_2014a,berg_2019a,kergosien_2023a}, where non-emergency patient transport and emergency transport are considered as an integrated optimization problem in order to utilize transport resources optimally.

\medskip

Figure~\ref{fig:goods} shows the distribution of the transported items. Due to the inclusion criteria of our review, all examined papers consider non-emergency patient transport, denoted by \emph{Patient}. Only a few papers~\cite{kallrath_2005a,melachrinoudis_2005a,vancroonenburg_2018a,Baermann2024} combine patient transport with the transport of material such as blood products, wheelchairs, or medication. Similarly, a few papers~\cite{kergosien_2014a,berg_2019a,kergosien_2023a} combine planned transports of non-emergency patients with transports of emergency patients by emergency medical services (EMS). This shows that models combining planned patient transport with other transport tasks in the hospital are rather uncommon.

\subsection{Transport fleet, depots, and constraints}
We now consider the characteristics of the transport fleet and the depots in the relevant literature as well as the considered constraints.

\begin{figure}[t]
  \begin{subfigure}[t]{\textwidth}
    \centering
    \includegraphics[width=\linewidth]{fleet_cons.png}
    \caption{Distribution of fleet characteristics ($\alpha$).}
        \label{fig:fleet_cons}
        \footnotesize
         \vspace{0.5em}
    \begin{tabular}{ll}
        No transport vehicles & \cite{nickel_2000a,fiegl_2009a,turan_2011a,schmid_2014a,elmbach_2015a,vancroonenburg_2018a,seguin_2019a,elmbach_2019a,xiao_2022a,maka_2022a,Baermann2024}  \\
        Heterogeneous fleet & \cite{kallrath_2005a,melachrinoudis_2005a,hanne_2009a,beaudry_2010a,kergosien_2011a,detti_2017a,berg_2019a} \\
        Homogeneous fleet & \cite{kergosien_2014a,zhang_2015a,nasir_2022a,kergosien_2023a} \\
    \end{tabular}
    
  \end{subfigure}
  \begin{subfigure}[t]{\textwidth}
    \centering
    \includegraphics[width=\linewidth]{depot_cons.png}
    \caption{Distribution of depot characteristics ($\beta$).}
        \label{fig:depot_cons}
        \footnotesize
         \vspace{0.5em}
    \begin{tabular}{ll}
        Multiple depots & \cite{kallrath_2005a,melachrinoudis_2005a,hanne_2009a,beaudry_2010a,kergosien_2011a,turan_2011a,kergosien_2014a,schmid_2014a,elmbach_2015a,zhang_2015a,detti_2017a,vancroonenburg_2018a,elmbach_2019a,berg_2019a,maka_2022a}  \\
        No depots & \cite{fiegl_2009a,seguin_2019a,xiao_2022a,Baermann2024} \\
        Single depot & \cite{nickel_2000a,nasir_2022a,kergosien_2023a} \\
    \end{tabular}
    
  \end{subfigure}
  \caption{Distribution of different fleet characteristics and depot characteristics in the examined papers.}
  \label{fig:contraints_one}
\end{figure}

\begin{figure}  
    \begin{subfigure}[t]{\textwidth}
    \centering
        \includegraphics[width=\linewidth]{cap_cons.png}
    \caption{Distribution of capacity constraints.}
    \footnotesize
     \vspace{0.5em}
    \begin{tabular}{ll}
        Unit capacity & \cite{nickel_2000a,kergosien_2011a,turan_2011a,kergosien_2014a,schmid_2014a,elmbach_2015a,vancroonenburg_2018a,elmbach_2019a,berg_2019a,xiao_2022a,maka_2022a,kergosien_2023a,Baermann2024}  \\
        Capacitated & \cite{kallrath_2005a,melachrinoudis_2005a,fiegl_2009a,hanne_2009a,beaudry_2010a,zhang_2015a,detti_2017a,nasir_2022a} \\
        No information & \cite{seguin_2019a} \\
    \end{tabular}
     \label{fig:cap_cons}
  \end{subfigure}
   \begin{subfigure}[t]{\textwidth}
    \centering
        \includegraphics[width=\linewidth]{timewindows.png}
    \caption{Distribution of time window constraints.}
    \footnotesize
     \vspace{0.5em}
    \begin{tabular}{ll}
        Mixed time windows (TW mix) & \cite{beaudry_2010a,kergosien_2014a,schmid_2014a,elmbach_2015a,detti_2017a,vancroonenburg_2018a,berg_2019a,Baermann2024}  \\
        Soft time windows (TW soft) & \cite{kallrath_2005a,melachrinoudis_2005a,fiegl_2009a,elmbach_2019a,xiao_2022a,nasir_2022a,kergosien_2023a} \\
         Hard time windows (TW hard) & \cite{kergosien_2011a,turan_2011a,zhang_2015a,maka_2022a} \\
        No information & \cite{nickel_2000a,hanne_2009a,seguin_2019a} \\
    \end{tabular}
    
     \label{fig:tw_cons}
  \end{subfigure}
  \hfill
  \begin{subfigure}[t]{\textwidth}
    \centering
    \includegraphics[width=\linewidth]{time_cons.png}
    \caption{Distribution of time-related constraints.}
     \label{fig:time_cons}
     \footnotesize
      \vspace{0.5em}
    \begin{tabular}{ll}
        Staff time & \cite{hanne_2009a,beaudry_2010a,kergosien_2011a,kergosien_2014a,zhang_2015a,vancroonenburg_2018a,seguin_2019a,berg_2019a,maka_2022a,kergosien_2023a,Baermann2024}  \\
        Loading time & \cite{hanne_2009a,beaudry_2010a,turan_2011a,zhang_2015a,detti_2017a,elmbach_2019a,nasir_2022a,maka_2022a,kergosien_2023a} \\
        Travel time & \cite{kallrath_2005a,hanne_2009a,beaudry_2010a,zhang_2015a,detti_2017a,maka_2022a} \\
    \end{tabular}
    
  \end{subfigure}
    \caption{Distribution of capacity constraints, time window constraints, and other time-related constraints in the examined papers.}
  \label{fig:contraints_two}
\end{figure}

\medskip

Figure~\ref{fig:fleet_cons} shows the distribution of fleet characteristics. Most papers do not consider transport with vehicles, i.e., they use a vehicle-free fleet. The reason for this is certainly the focus of our literature search on intra-hospital patient transport. Indeed, the papers that use a vehicle-free fleet are a subset of the papers categorized as inside hospitals in Figure~\ref{fig:locations}. The papers~\cite{hanne_2009a,kallrath_2005a,beaudry_2010a}, which only mention transport inside hospitals and use vehicles for transport, consider large hospital complexes -- sometimes on mile-long road networks~\cite{beaudry_2010a}. 
In most papers considering vehicles, ambulances are used for transport. However, ambulances for non-emergency transport (so-called Basic Life Support and Advanced Life Support Ambulances~\cite{melachrinoudis_2005a}) are often mentioned, which leads to an increased use of heterogeneous fleets, alongside fleets that also differ in their capacity and tasks. In addition, some studies~\cite{kergosien_2011a,beaudry_2010a,hanne_2009a,melachrinoudis_2005a} model a heterogeneous fleet to reflect real-world conditions.

\medskip

With respect to depots, Figure~\ref{fig:depot_cons} shows a clear dominance of multi-depot problems. It should be noted, however, that a large proportion of papers do not give specific reasons for using multi-depot approaches. Some papers~\cite{melachrinoudis_2005a,beaudry_2010a} state that the use of a multi-depot approach is based on real data. The few papers using single-depot approaches use a central assembly point for the transporters or transport vehicles, which can also be seen as a result of real data~\cite{nickel_2000a}. For the problems that do not use depots, this modeling choice results from the specific approaches. For instance, Séguin et al.~\cite{seguin_2019a} assign responsibilities to the transporters for transport corridors between building sections, and in~\cite{xiao_2022a}, nurses begin their transport tasks from their assigned nursing ward, accompany patients to their appointments or pick them up from appointments, and then return directly to their ward. 

\medskip

During our interviews with hospital practitioners, we were told that the transporters do not have real depots, i.e., a room, but simply return to various places in the hospital complex where they can sit and wait between transports. This indicates that the establishment of assembly points or depots in a hospital is not exactly planned, but that tactically useful and central waiting points with seating facilities, for example, are established as depots by the transporters. The use of multiple depots can also be advantageous in terms of minimizing travel distances and waiting times and, thus, lead to overall gains in efficiency of the patient transport process.

\medskip
 Figure~\ref{fig:cap_cons} shows the distribution of capacity constraints. The majority of the papers assume a unit capacity, which is mainly due to the dominance of vehicle-free transport. In addition, ambulances, which in some papers~\cite{berg_2019a,kergosien_2014a,kergosien_2023a} also provide non-emergency patient transport in addition to emergency transport, are often assumed to have unit capacity. Capacitated problems are usually associated with a heterogeneous fleet, but vehicle-free problems are also modeled as capacitated problems. This is due to the fact that transport personnel can usually only push one bed or wheelchair at a time, but can accompany two walking patients simultaneously~\cite{fiegl_2009a}.

\medskip

Even though patients often require special equipment such as oxygen supply or monitoring equipment during transports, such additional equipment requirements are only considered in four of the 22 examined papers~\cite{kallrath_2005a,hanne_2009a,beaudry_2010a,kergosien_2011a}. 
Different priorities for transport requests are included more frequently, with nine papers considering them~\cite{nickel_2000a,fiegl_2009a,beaudry_2010a,kergosien_2011a,turan_2011a,kergosien_2014a,elmbach_2015a,berg_2019a,Baermann2024}, whereas different skill levels of transport personnel are only considered in six cases~\cite{beaudry_2010a,zhang_2015a,vancroonenburg_2018a,berg_2019a,kergosien_2023a,Baermann2024}. The latter, in particular, seems surprising since skill level requirements such as the need for a patient under the influence of medication to be transported by medically trained personnel are not only important for patient safety, but may also constitute a legal requirement. Similarly, cases where patients require isolated transport (e.g., due to infectious diseases) are also considered in only six of the 22 examined papers~\cite{kallrath_2005a,hanne_2009a,beaudry_2010a,kergosien_2011a,detti_2017a,Baermann2024}, even though this requirement may also have important safety implications as well as logistical implications such as specialized protective clothing that needs to be provided for the transport personnel performing these transports. Moreover, isolated transports may also imply the need for disinfection or special cleaning of vehicles after completion of a transport, which can extend the unavailability of vehicles.

\medskip
Figure~\ref{fig:tw_cons} illustrates how frequently the different types of time window constraints are used in the examined papers. As the figure shows, 15 out of the 19 papers for which information about time window constraints is available use either mixed or soft time windows, while hard time windows are only used in four papers. This can largely be attributed to the high frequency of delays and unforeseen events in patient transport, which mean that problems with hard time windows might easily become infeasible in high utilization scenarios. In problems with mixed time window constraints, soft pick-up times are often implemented to accommodate potential delays due to preceding appointments, while delivery times are more frequently set as hard constraints to ensure timely continuation of subsequent medical procedures and to prevent equipment and personnel at the destination from becoming idle due to late arrivals of patients. Alternatively, some papers assume that the patient is ready for pick-up at a fixed time, while allowing flexibility in the delivery time. In cases where both pick-up and delivery times are flexible, penalties are commonly applied for early arrival at the pick-up location or late arrival at the delivery location.

\medskip

\enlargethispage{0.75\baselineskip}

Figure~\ref{fig:time_cons} considers other time-related constraints, where staff time constraints, i.e., constraints on working and break times of transport staff, occur most frequently, followed by loading time and travel time constraints. It should be noted, however, that these constraints are not mutually exclusive, i.e., staff time constraints, loading time constraints, and travel time constraints may all be considered together in one problem, and they might also be combined with time window constraints. 
Loading times appear in the examined papers as times for the disinfection of ambulances~\cite{beaudry_2010a,kallrath_2005a} or as times for preparation and follow-up work~\cite{turan_2011a,zhang_2015a,detti_2017a}. For problems with capacities strictly larger than the unit capacity, where several patients are transported and a loaded patient is not necessarily at their destination at the next stop, travel times are more important. For example, in the case of combined transports of several patients, the maximum travel time is sometimes limited by a multiple of the direct travel time~\cite{zhang_2015a}.

\subsection{Objectives}

Figure~\ref{fig:objective_characteristics} shows the distribution of optimization objectives. Here, as already mentioned, cost minimization may comprise many different kinds of costs, some of which arise by assigning costs to constraint violations. The most important costs considered in the relevant literature are the total costs of transport~\cite{kergosien_2023a, zhang_2015a}, the costs of vehicles~\cite{nasir_2022a}, and personnel costs~\cite{seguin_2019a}. 
Concerning the optimization of waiting times, some papers minimize waiting times for patients~\cite{kallrath_2005a,melachrinoudis_2005a,hanne_2009a,beaudry_2010a,turan_2011a,kergosien_2014a,schmid_2014a,zhang_2015a,detti_2017a,xiao_2022a,nasir_2022a,maka_2022a}, while others minimize vehicle waiting times~\cite{kallrath_2005a,hanne_2009a,beaudry_2010a,turan_2011a,schmid_2014a,zhang_2015a,detti_2017a,nasir_2022a}. 

\begin{figure}
    \centering
   \includegraphics[width=\textwidth]{objective.png}
    \caption{Distribution of different objectives ($\delta$) considered in the examined papers ordered as in the five-field notation (WT = waiting time, TT = travel time).}
    \label{fig:objective_characteristics}
    \footnotesize
     \vspace{0.5em}
    \begin{tabular}{ll}
         Cost &  \cite{nickel_2000a,kallrath_2005a,melachrinoudis_2005a,hanne_2009a,beaudry_2010a,kergosien_2011a,turan_2011a,kergosien_2014a,schmid_2014a,zhang_2015a,detti_2017a,seguin_2019a,nasir_2022a,maka_2022a,kergosien_2023a} \\
         WT patient &  \cite{kallrath_2005a,melachrinoudis_2005a,hanne_2009a,beaudry_2010a,turan_2011a,kergosien_2014a,schmid_2014a,zhang_2015a,detti_2017a,xiao_2022a,nasir_2022a,maka_2022a} \\
         WT vehicle &  \cite{kallrath_2005a,hanne_2009a,beaudry_2010a,turan_2011a,schmid_2014a,zhang_2015a,detti_2017a,nasir_2022a} \\
         TT vehicle &  \cite{kallrath_2005a,hanne_2009a,beaudry_2010a,turan_2011a,zhang_2015a,detti_2017a,vancroonenburg_2018a,berg_2019a,nasir_2022a,maka_2022a,kergosien_2023a} \\
         TT patient &  \cite{kallrath_2005a,melachrinoudis_2005a,hanne_2009a,beaudry_2010a,zhang_2015a,detti_2017a,xiao_2022a,nasir_2022a,maka_2022a} \\
         TT empty &  \cite{nickel_2000a,beaudry_2010a,turan_2011a,kergosien_2014a,schmid_2014a,detti_2017a,nasir_2022a,Baermann2024} \\
         Ergonomic burden &  \cite{elmbach_2015a,elmbach_2019a} \\ 
         Flow time &  \cite{fiegl_2009a} \\
    \end{tabular}
    
\end{figure} 

\medskip
It can be seen that there are differences between waiting times and travel times with respect to the prioritization of patients and vehicles.
For example, all of the 12 papers considering waiting times include waiting times for patients. A potential reason is that personnel are required to monitor patients during transports, so minimizing patient waiting times not only improves the quality of service but also saves human resources. Waiting times for vehicles also play an important role, but are only considered in eight out of 12 papers modeling waiting times~\cite{kallrath_2005a,hanne_2009a,beaudry_2010a,turan_2011a,schmid_2014a,zhang_2015a,detti_2017a,nasir_2022a}. Their consideration is sometimes explained by the low availability of vehicles~\cite{beaudry_2010a,turan_2011a,detti_2017a,zhang_2015a}.

\medskip

When optimizing travel times, a different approach can be seen in the prioritization of patients and vehicles. In 11 out of 17 cases~\cite{kallrath_2005a,hanne_2009a,beaudry_2010a,turan_2011a,zhang_2015a,detti_2017a,vancroonenburg_2018a,berg_2019a,nasir_2022a,maka_2022a,kergosien_2023a}, the focus of travel time optimization is on vehicles, which is again motivated by the low availability of vehicles or the goal of using vehicles as efficiently as possible. This is underpinned by the fact that eight out of the 17 papers minimizing travel times also try to minimize empty travel times of vehicles~\cite{nickel_2000a,beaudry_2010a,turan_2011a,kergosien_2014a,schmid_2014a,detti_2017a,nasir_2022a,Baermann2024}. The minimization of travel times for patients, which is considered in nine out of the 17 papers~\cite{kallrath_2005a,melachrinoudis_2005a,hanne_2009a,nasir_2022a,xiao_2022a,detti_2017a,zhang_2015a,beaudry_2010a,maka_2022a}, is slightly less common in the examined papers with a travel time objective than the minimization of vehicle travel times. Reasons for this may be that travel times for transports within hospitals are often considered irrelevant~\cite{fiegl_2009a}, as well as that many of the studied papers work with unit capacities, where a direct trip from origin to destination can be assumed and, thus, travel times for patients are fixed unless uncertain travel times are considered.

\medskip

Concerning the simultaneous consideration of multiple objectives, 12 out of 22 papers incorporate multiple objectives by using a weighted-sum objective function~\cite{kallrath_2005a, melachrinoudis_2005a, hanne_2009a, beaudry_2010a, turan_2011a, schmid_2014a, zhang_2015a, detti_2017a, xiao_2022a, nasir_2022a, maka_2022a, kergosien_2023a}. Only two papers seem to adopt real multi-objective approaches, namely~\cite{nickel_2000a} by minimizing several objectives including the number of required transport staff and the time the transport staff spend walking without a patient, and~\cite{kergosien_2014a} by using a lexicographic approach for handling multiple objectives such as transport delays, staff overtime, and traveled distance. In contrast, eight out of 22 studies only use a single objective function~\cite{fiegl_2009a, kergosien_2011a, elmbach_2015a, elmbach_2019a, vancroonenburg_2018a, seguin_2019a, Baermann2024, berg_2019a}.  

\subsection{Uncertainty and uncertainty handling}
A wide range of uncertainties such as uncertain requests, uncertain travel times, uncertain availability of staff or vehicles, and many more are conceivable in patient transport problems. The importance of taking uncertainty into account is highlighted in some papers~\cite{melachrinoudis_2005a,schmid_2014a,elmbach_2015a}, which describe that their models have not made it into practical use due to a lack of uncertainty handling (the relationship between the inclusion of uncertainty and the level of practical implementation is considered in more detail in Section~\ref{sec:gap}). It is noteworthy, however, that, as illustrated in Figure~\ref{fig:uncertainties}, half of the examined papers do not consider any uncertainty at all. As the figure shows, the most frequently considered kind of uncertainty is uncertain transport requests, which means that not all transport requests are known in advance, but new requests may arrive dynamically over the considered time horizon. 
The paper by Turan et al.~\cite{turan_2011a} states that a major reason for the increasing consideration of uncertain requests is the observation that a large proportion of transport requests are not communicated to the transport service as soon as they become known. Avoiding this delay could potentially allow more time for incorporating arriving requests in the existing transport schedule. 
As can also be seen from Figure~\ref{fig:uncertainties}, uncertain travel times are considered much less frequently than uncertain requests. A potential reason could be that the focus of the examined papers is mostly on transport within the hospital and, as stated in~\cite{fiegl_2009a}, routes within a hospital are often very short and therefore do not play a relevant role. Overall, the consideration of uncertainty means that problems that are already hard to solve become even more complex.

\begin{figure}
  \begin{subfigure}[t]{\textwidth}
    \centering
    \includegraphics[width=\linewidth]{uncertain.png}
    \caption{Distribution of different kinds of uncertainty}.
        \label{fig:uncertainties}
         \vspace{0.5em} 
         \footnotesize
    \begin{tabular}{ll}
         No uncertainty &  \cite{melachrinoudis_2005a,turan_2011a,schmid_2014a,elmbach_2015a,elmbach_2019a,zhang_2015a,detti_2017a,seguin_2019a,xiao_2022a,maka_2022a,nasir_2022a}  \\
         Uncertain requests &  \cite{nickel_2000a,kallrath_2005a,fiegl_2009a,hanne_2009a,beaudry_2010a,kergosien_2011a,kergosien_2014a,vancroonenburg_2018a,berg_2019a,kergosien_2023a,Baermann2024} \\
         Uncertain travel times &  \cite{nickel_2000a,hanne_2009a,beaudry_2010a,berg_2019a} \\
    \end{tabular}
    
  \end{subfigure}
  \hfill
  \begin{subfigure}[t]{\textwidth}
    \centering
    \includegraphics[width=\linewidth]{hand_uncertain.png}
    \caption{Distribution of different methods used to deal with uncertainty.}
     \label{fig:handling-strategies}
     \footnotesize
      \vspace{0.5em}
    \begin{tabular}{ll}
        Re-optimization & \cite{kallrath_2005a,hanne_2009a,beaudry_2010a,kergosien_2011a,kergosien_2014a,vancroonenburg_2018a,berg_2019a,kergosien_2023a,Baermann2024} \\
        Online optimization & \cite{nickel_2000a,kallrath_2005a,fiegl_2009a} \\
        Stochastic optimization & \cite{nickel_2000a} \\
    \end{tabular}
    
  \end{subfigure}
  \caption{Distribution of different kinds of uncertainty considered and uncertainty handling methods used in the examined papers.}
  \label{fig:unc-unc-hand}
\end{figure}

\medskip

Various approaches to handle the mentioned types of uncertainty are used in the examined literature. Figure~\ref{fig:handling-strategies} shows that re-optimization is the most frequently used method for dealing with uncertainty. The only exceptions are Nickel and Tenfelde~\cite{nickel_2000a}, who use a combination of stochastic and online optimization,  and Kallrath~\cite{kallrath_2005a} and Fiegl and Pontow~\cite{fiegl_2009a}, who also use online optimization approaches. In re-optimization, either the optimization process is re-started from scratch in order to find a new solution, or parts of the existing solution are re-optimized as soon as problem parameters, such as the set of requests, change.
One reason for the frequent use of re-optimization as a strategy for dealing with uncertainty is that re-optimization can react flexibly to a constantly changing number of requests. Moreover, as in~\cite{beaudry_2010a,hanne_2009a,vancroonenburg_2018a}, re-optimization allows to integrate new requests quickly into existing solutions. 

\subsection{Modeling and solution methods}
Figure~\ref{fig:sol_methods} shows the distribution of different modeling and solution methods used in the examined papers. Note that these are not mutually exclusive, so it is possible that several methods are used in a single paper either independently or in combination. Most frequently, (mixed-) integer programming is used to model patient transport problems in the relevant papers. One advantage of (mixed-) integer programming for modeling is that it has interdisciplinary recognition and allows to represent complex details of an optimization problem in an exact way~\cite{floudas_2005a}. Other approaches include modeling a patient transport problem as a quadratic assignment problem~\cite{nickel_2000a} or as a job scheduling problem~\cite{elmbach_2015a,fiegl_2009a}.

\begin{figure}
    \centering
   \includegraphics[width=\textwidth]{solmeth.png}
    \caption{Distribution of different modeling and solution methods used in the examined papers. The term \emph{Others} encompasses all solution methods that do not fit into the previous categories.}
    \label{fig:sol_methods}
     \vspace{0.5em}
     \footnotesize
    \begin{tabular}{ll}
        Mixed-integer programming (MIP) & \cite{kallrath_2005a,melachrinoudis_2005a,fiegl_2009a,turan_2011a,schmid_2014a,elmbach_2015a,zhang_2015a,detti_2017a,seguin_2019a,elmbach_2019a,berg_2019a,xiao_2022a,nasir_2022a,maka_2022a,kergosien_2023a,Baermann2024}  \\
        Improvement heuristics & \cite{melachrinoudis_2005a,beaudry_2010a,kergosien_2011a,kergosien_2014a,schmid_2014a,elmbach_2015a,zhang_2015a,vancroonenburg_2018a,elmbach_2019a} \\
        Greedy \& insertion heuristics & \cite{nickel_2000a,hanne_2009a,beaudry_2010a,schmid_2014a,detti_2017a,vancroonenburg_2018a,xiao_2022a,nasir_2022a} \\
        Metaheuristics & \cite{kallrath_2005a,hanne_2009a,zhang_2015a,xiao_2022a} \\
        Problem-specific heuristics & \cite{kallrath_2005a,fiegl_2009a,kergosien_2023a} \\
        Others & \cite{nickel_2000a,kergosien_2011a,nasir_2022a} \\
    \end{tabular}
    
\end{figure} 
\medskip

Concerning solution methods, the examined papers often combine other methods (e.g., (mixed-) integer programming) with heuristics, especially with improvement heuristics~\cite{melachrinoudis_2005a,beaudry_2010a,kergosien_2011a,kergosien_2014a,schmid_2014a,elmbach_2015a,zhang_2015a,vancroonenburg_2018a,elmbach_2019a} or greedy and insertion heuristics~\cite{nickel_2000a,hanne_2009a,beaudry_2010a,schmid_2014a,detti_2017a,vancroonenburg_2018a,xiao_2022a,nasir_2022a}. Here, greedy heuristics are often used to generate high-quality initial solutions, which can reduce the computation time required by other solution methods such as (mixed-) integer programming solvers. Other papers use insertion heuristics to solve the considered problem or to find partial solutions~\cite{zhang_2015a,detti_2017a,vancroonenburg_2018a,hanne_2009a,beaudry_2010a}. Moreover, insertion heuristics are also used for re-optimization because they allow a quick integration of new requests into existing transport schedules. Improvement heuristics are frequently used in the examined papers to improve initial solutions generated by greedy heuristics or time-limited runs of (mixed-) integer programming solvers, which can reduce the optimization time significantly. Other notable approaches used include metaheuristics~\cite{kallrath_2005a,hanne_2009a,zhang_2015a,xiao_2022a} and problem-specific heuristics~\cite{kallrath_2005a,fiegl_2009a,kergosien_2023a}. 

\subsection{Data and level of real-world implementation}

\begin{figure}
  \begin{subfigure}[t]{\textwidth}
    \centering
    \includegraphics[width=\linewidth]{realdata.png}
    \caption{Distribution of the use of real data for model development.}
    \label{fig:data-source}
    \vspace{0.5em}
    \footnotesize
    \begin{tabular}{ll}
        Full real data & \cite{kallrath_2005a,melachrinoudis_2005a,fiegl_2009a,hanne_2009a,beaudry_2010a,kergosien_2011a,zhang_2015a,detti_2017a,seguin_2019a,berg_2019a,xiao_2022a,Baermann2024}  \\
        No real data & \cite{nickel_2000a,turan_2011a,kergosien_2014a,schmid_2014a,elmbach_2015a,elmbach_2019a,maka_2022a,kergosien_2023a} \\
        Real hospital layout / location data & \cite{vancroonenburg_2018a,nasir_2022a} \\
    \end{tabular}
    
  \end{subfigure}
  \hfill
  \begin{subfigure}[t]{\textwidth}
    \centering
    \includegraphics[width=\linewidth]{practise.png}
    \caption{Distribution of the level of practical implementation achieved.}
    \label{fig:integration-level}
    \vspace{0.5em}
    \footnotesize
    \begin{tabular}{ll}
        Testing with full real data & \cite{kallrath_2005a,melachrinoudis_2005a,fiegl_2009a,beaudry_2010a,kergosien_2011a,zhang_2015a,detti_2017a,seguin_2019a,berg_2019a,xiao_2022a,Baermann2024}  \\
        No testing with full real data & \cite{nickel_2000a,turan_2011a,kergosien_2014a,schmid_2014a,elmbach_2015a,elmbach_2019a,vancroonenburg_2018a,nasir_2022a,maka_2022a,kergosien_2023a} \\
        Used in practice & \cite{hanne_2009a} \\
    \end{tabular}
    
  \end{subfigure}
  \caption{Distribution of the use of real data and the level of practical implementation achieved in the examined papers.}
  \label{fig:real-data}
\end{figure}

A real-world patient transport problem instance consists of various data. First of all, the pick-up and delivery locations and, if applicable, the locations of the depots should be specified. This can be done using coordinates or a distance matrix. For problems that include capacities or a heterogeneous fleet, the vehicle or transporter properties must also be specified, as well as the starting depots in the case of multi-depot problems. Finally, the data for the requests must be provided. In the examined literature, a distinction is only made between real location data and real requests. Figure~\ref{fig:data-source} shows that more than half of the proposed models have been developed using full real data. In contrast, Vancroonenburg et al.~\cite{vancroonenburg_2018a} and Nasir et al.~\cite{nasir_2022a} only use real hospital layouts or real location data. The reason for this is that layouts are usually available for the hospital buildings, or distances outside hospitals can be easily determined using Geographic Information Systems (GIS)~\cite{berg_2019a}. None of the examined papers use real requests as the only piece of real data. The reason is that there are often no records of requests~\cite{nickel_2000a} and that the origin and the destination are part of the information that a request must contain, so using real request data for model development usually implies the use of a real hospital layout.

\medskip

It is also noticeable that, apart from~\cite{vancroonenburg_2018a}, none of the examined papers provide the considered problem instances. Since the paper~\cite{vancroonenburg_2018a} works with real hospital data but artificially generated requests, this means that no benchmark instances consisting of completely real data are currently available for future research.

\medskip

Moreover, as Figure~\ref{fig:integration-level} shows, except for~\cite{hanne_2009a}, none of the developed models have been implemented in practice. The reasons for this stated in the papers range from performance issues over missing practice-relevant constraints to uncertainty about the cost-benefit aspect of integration in practice. A more detailed overview and categorization of the mentioned reasons for missing practical implementation is provided in Section~\ref{sec:gap}.
\subsection{Cross comparisons}\label{sec:cross}

We now investigate connections that result from cross comparisons between the different aspects analyzed in the previous subsections.

\medskip

First, an interesting relationship between the modeling of transport capacity (see also Table~\ref{tab:paper_five_field}) and maximum travel time constraints can be observed. Of the 13 papers using unit capacity, only~\cite{maka_2022a} includes maximum travel time constraints. In contrast, five of the eight papers using larger capacities do so. This is because unit capacity problems assume direct transport between origin and destination, which means optimizing travel times for patients is irrelevant. Moreover, nine out of 13 papers studying unit capacity problems consider transport inside hospitals, where distances and travel times are so marginal that they are of little relevance~\cite{fiegl_2009a}. When larger capacities are considered, however, longer tours involving simultaneous transport of several patients can be formed, which makes limiting patient travel times much more relevant. 

\medskip

A further relevant cross-comparison shows that papers considering uncertainty are more likely to include constraints on staff working hours and breaks (eight of 11 papers). This could be explained by the higher risk of working or break time violations that results from uncertain travel times or requests, which makes it more important to explicitly incorporate constraints that avoid these violations into the problems.

\medskip

Lastly, we examine patient transport types based on origin-destination pairs. It turns out that, with the exception of~\cite{melachrinoudis_2005a}, none of the papers consider all three types of non-emergency patient transport (transports inside hospitals, transports between hospitals (transfers), and transports from home to hospital or hospital to home (discharges)). In contrast, several of the examined papers consider transfers and discharges of patients together~\cite{berg_2019a,detti_2017a,kergosien_2014a,nasir_2022a,kergosien_2023a}, in some cases also in combination with the transport of patients from accident locations to the hospital~\cite{kergosien_2014a,berg_2019a,kergosien_2023a}. Finally, it is also noticeable that there is a connection between the transports considered and the use of real data. The papers dealing exclusively with patient transports within hospitals~\cite{nickel_2000a,turan_2011a,kallrath_2005a,beaudry_2010a,fiegl_2009a,schmid_2014a,xiao_2022a,elmbach_2015a,elmbach_2019a,maka_2022a,Baermann2024,vancroonenburg_2018a,seguin_2019a} are less likely to use real data. One reason for this could be the practice of registering transports by telephone in some hospitals, which does not include data storage due to the lack of digitization~\cite{hanne_2009a,nickel_2000a}.
\section{Research gaps and agenda}\label{sec:gap}
This section identifies research gaps resulting from our literature analysis and derives an agenda for future research. In order to optimally align our findings with requirements from hospital practice, we conducted semi-structured interviews with people responsible for the organization of patient transport in five selected hospitals of different sizes in Germany and the Netherlands. Table~\ref{tab:hosp} below summarizes the bed capacity, total staff, and the number of full-time equivalent transport staff in each hospital. The English translation of our interview guide is provided in~\ref{sec:appa}. Besides specific questions that resulted from our literature analysis, we also asked open questions regarding the general organization of patient transport, challenges encountered, and possible potential for improvement, so that our interview partners had the opportunity to also mention aspects not considered by us or in the existing literature.

\begin{table}[h!]
    \centering
    \begin{tabular}{|c|c|c|c|}
        \hline
        \textbf{Hospital} & \textbf{Beds} & \textbf{Staff} & \textbf{Transport staff} \\
        \hline
        Hospital 1 & 1700 & 16000 & 36 \\
        Hospital 2 & 1160 & 6600 & 22 \\
        Hospital 3 & 660 & 1900 & 15 \\
        Hospital 4 & 560 & 1500 & 5 \\
        Hospital 5 & 475 & 1800 & 14 \\
        \hline
    \end{tabular}
    \caption{Characteristics of the selected hospitals.}\label{tab:hosp}
\end{table}

\medskip

In the following, we first provide a short description of the organization of patient transport that is currently in place at the selected hospitals (Subsection~\ref{sec:gap_current_situation}). The following subsections then compare the research gaps and results of our literature analysis with the statements made during the interviews in order to obtain a comprehensive research agenda. In particular, we identify important aspects of the practical problem that are currently missing from the patient transport problems studied in the literature (Subsection~\ref{sec:gap_missing}), research gaps related to the integration of patient transport with other planning problems in hospitals (Subsection~\ref{sec:gap_integration}), problems related to data availability (Subsection~\ref{sec:gap_data}), and possible reasons for the low number of successful practical implementations of existing patient transport models and algorithms (Subsection~\ref{sec:gap_practical}).

\subsection{Current situation at the selected hospitals}\label{sec:gap_current_situation}
In one of the hospitals, patient transport is still organized without a specialized software solution. Instead, a roster is manually created every morning, assigning each transporter to a specific kind of transport requests, such as handling all transports to and from the radiology department. This roster, along with the transporters' phone numbers, is then distributed to the relevant hospital departments. To request a transport, the initiating person (e.g., a nurse on a nursing ward) calls the transporter responsible for this specific kind of transport request and registers the transport by phone. Transports performed by medical staff (e.g., nurses) are initiated by phone calls between the pick-up and delivery locations. During these calls, arrangements are made regarding who will transport the patient, whether the patient will be picked up and brought to the destination, and about any specific requirements or instructions for the transport, such as medical considerations or timing.

\medskip

Four out of five of the hospitals, however, use dedicated software that is directly connected to the hospital information system in order to coordinate patient transports. In one case, in addition to patient transport, the software also coordinates the transport of materials, medicines, and blood samples. In order to initiate a patient transport to be performed by a transporter via the software, the department or nursing ward from which the transport originates submits a transport request in the software, which is then automatically forwarded by the system to one or several transporters. The transporters carry a smartphone with them on which a dedicated app informs them about incoming transport requests, which they can then accept directly via the app. In the case of one hospital, to carry out a transport, a transporter has to scan QR codes at the pick-up location and on the patient's wristband when the transport starts, and afterwards at the drop-off location and again on the patient's wristband in order to allow the software to register the pick-up and delivery times and keep track of the patient's location. In the case of another hospital, the patient is identified by a five-point system at the start of a transport, which includes checking the patient's location on the nursing ward layout and on the room allocation plan, as well as checking the patient's number, name, and date of birth. After the transport, the successful delivery is then confirmed via the app. However, neither the software nor the smartphone app include any assistance for the routing of the transporters. The decision about which transporters are notified of an incoming transport request is managed automatically by the software, based on a proprietary, black-box algorithm. According to the available information, the assignments made by the algorithm are chosen so as to balance the number of requests across transporters during periods of low transport demand, but the algorithm does not perform any optimization beyond that. Moreover, the software systems only coordinate transports carried out by patient transporters, while transports to be performed by medical personnel, such as nurses, are handled manually, similar to the phone call-based process described above.

\subsection{Missing aspects in patient transport}\label{sec:gap_missing}
This subsection highlights important aspects of patient transport pointed out by our interview partners that are currently missing from or not sufficiently considered in the problems studied in the literature and / or the practical processes currently used for organizing patient transport.

\medskip

One important aspect mentioned by all of our interview partners is considering the workload and ergonomic burden of the transport staff -- in particular, the dedicated patient transporters. Due to personnel shortages, minimizing the workload of these employees and distributing it fairly is very important according to all of our interview partners. This aspect, however, is not sufficiently considered so far either in practice or in the literature.

\medskip

For example, the hospital that still coordinates patient transport without a specialized software solution states that their current system does not allow for any systematic consideration or balancing of individual workloads of transporters. This is despite the fact that the employees of the transport service often show a large number of physical symptoms resulting from excessive workload and ergonomic burden, which often leads to a large number of sick days and, in turn, to an even larger workload for the remaining staff.
To tackle this problem, both better planning that takes individual workloads into account and the use of pushing aids for patient transport are mentioned as very desirable by four out of five interview partners. 
One of the hospitals already describes using a complete fleet of pushing aids (so-called \emph{bed movers}). These devices make transport more comfortable for the transporters and reduce their ergonomic burden because the transporters can travel on them. In addition, this hospital uses an interesting workload-balancing strategy based on rotating responsibilities. Here, the transport staff belong to a logistics department that is responsible for material transport, patient transport, as well as bed cleaning.
The work in bed cleaning is described as being less physically demanding than patient transport, even though pushing aids are used. During busy periods, all logistics department staff work in patient transport; during periods of low workload, the employees switch to bed cleaning according to a predefined roster. The roster rotates so that each employee is assigned to a different responsibility each day to balance the workload. Another interview partner from a hospital that uses software to coordinate patient transports states that the software already balances the number of transports to be performed among the transporters during periods of low demand. However, neither of these two approaches considers the specific workload and ergonomic strain resulting from individual transport requests. Apart from these basic balancing methods, no approaches for minimizing the workload (e.g., by suggesting transport tours that minimize transport distances or times) are currently in place.

\medskip

We believe that this feedback from practice demonstrates the high relevance of including transport staff workload and ergonomic burden in patient transport problems. In the patient transport literature, however, this is currently only considered by von Elmbach et al.~\cite{elmbach_2015a,elmbach_2019a}. They investigate factors that lead to ergonomic stress for transporters to minimize the maximum ergonomic burden. In doing so, they look at an increasingly important aspect for the sustainable use of transport and nursing staff in hospitals, who face increasing workloads due to staff shortages and an increasing demand for health care. Future research should therefore put a larger emphasis on both ergonomic aspects of patient transport and on workload distribution among the transport staff, which could also contribute to breaking the downward spiral between staff shortages and overloading of available employees.

\medskip

Another essential aspect of patient transport described by all our interview partners is the qualification levels of transport staff. Some transports, such as those where patients are under the influence of medication, should be performed by medically trained personnel (e.g., trips to and from the operating room), while other requests (e.g., transports to radiology and other diagnostic units) can be handled by medically untrained transporters. 
For transports to and from the operating room, for example, different qualification levels are often taken into account in practice by administering medication for surgeries on site, so that patients can be transported to the operating room by medically untrained staff. After surgery, all patients are divided into low- and high-risk patients, so that untrained transporters can transport low-risk patients, while all other transports must be carried out by medically trained personnel from the destination nursing wards. These two kinds of transports, however, are not yet planned jointly in a centralized manner since, even in the hospital that already uses a specialized software solution to coordinate patient transports performed by medically untrained transporters, all transports performed by medical personnel are still planned separately without any software support. In the literature, there also seems to be a large gap with respect to the consideration of skill levels of transport staff, which are so far only considered in~\cite{beaudry_2010a,zhang_2015a,vancroonenburg_2018a,berg_2019a,kergosien_2023a,Baermann2024}. Therefore, including skill levels and the joint planning of medically trained and untrained transport staff in patient transport problems represents an important topic for future research with large practical relevance. 

\medskip
A further relevant but not yet sufficiently considered aspect of patient transport is to integrate transports of material such as bed sheets, blood samples, medication, or bandages. While all of our interview partners acknowledge the potential of joint patient and material transports by patient transporters, this is only very rarely implemented in the selected hospitals so far. One of the hospitals states that the integration of material transports into patient transports is not currently practiced and is not planned for the future. This has mainly administrative reasons resulting from different cost centers for patient transport and material transport that require a strict separation of the two kinds of transports. While another hospital already has a joint logistics department that handles patient transport as well as material transport, they nevertheless perform the two types of transports separately. Other hospitals confirm that a few transports of materials are already carried out by untrained transporters mainly responsible for patient transport and that further material transport tasks, such as transports of surgical instruments to and from the sterilization department, are planned to be integrated with patient transports in the future. The challenge here is the need for additional capacity of the patient transport service.
Moreover, one of the hospitals states that transports of wheelchairs back to their desired storage locations are not yet considered during patient transport planning, which results in a lack of wheelchairs when picking up patients, or in the parking of wheelchairs in inconvenient places. 

\medskip

In summary, the integration of material transport and the transport of equipment such as wheelchairs is an important practical aspect that offers large potential for better staff utilization and cost savings. In the patient transport literature, however, only a few of the examined papers consider this aspect~\cite{kallrath_2005a,melachrinoudis_2005a,vancroonenburg_2018a,Baermann2024}, and no results are available on the practical implementation of models for integrated patient and material transport.

\subsection{Integration with other planning problems}\label{sec:gap_integration}
We now consider the integration of patient transport with other planning problems in hospitals. It is well-known that patients must be present for a large part of the activities performed in a hospital, and that delayed arrivals of patients for examination, surgery, or other procedures can lead to significant disruptions in day-to-day operations. It is, therefore, surprising that, as also observed in~\cite{rachuba_2023a}, there are only few papers that deal with the integration of patient transport with other planning problems in hospitals.

\medskip

In our literature search, we found only two papers that consider this integration explicitly~\cite{nickel_2000a,schmid_2014a}. The integration into hospital layout planning considered in~\cite{nickel_2000a}, however, cannot be assigned to the operational planning and optimization of daily hospital activities. Thus, the paper by Schmid and Doerner~\cite{schmid_2014a} is the only paper in our literature search that integrates patient transport with other operational planning problems in hospitals -- in this case, with operating room scheduling. The observation that patient transport is first integrated with operating room scheduling is not surprising since a large part of integrated planning in hospitals is aimed at optimizing the operating room, which represents the largest cost and revenue driver in hospitals~\cite{rachuba_2023a}. However, since patient transport interacts with a large variety of planning problems in hospitals, it should ideally be integrated with all of these problems arising in many different departments. One example of how patient transport and appointment scheduling in diagnostic units such as radiology could be integrated is the generation of a patient transport request as soon as the diagnostic procedure is scheduled. This information would then be incorporated into the determination of the required number and type of transport personnel. Overall, we conclude that, in addition to further research on patient transport itself, integration with other operational planning problems in hospitals also represents an important area for future research. 

\medskip

Moreover, as already remarked in Section~\ref{sec:cross}, there is no work that examines patient transport in hospitals, patient transport between hospitals, and the transport of patients after discharge jointly in one model. As these different kinds of patient transports clearly interact and sometimes involve the same staff, this represents another research gap and also a practical challenge. For instance, integration of intra-hospital transports and transports after discharge could involve hospital staff taking a discharged patient to a pick-up area near an ambulance bay shortly before the ambulance arrives, thus saving the ambulance crew a significant amount of time, as they would not have to search for and pick up the patient from a nursing ward. One of the selected hospitals already using specialized software for patient transport planning explains that some external patient transports could be integrated into the software tool they use. Still, technical and management interfaces to external partners such as ambulance services and other hospitals are missing. 
Another hospital states that patients are not usually discharged directly from the nursing wards, but that most patients are transferred to a dedicated rehabilitation ward within the hospital complex after their stationary treatment. The discharge is then coordinated on the rehabilitation ward, which again means that no integration of discharges and patient transports within the hospital takes place. Therefore, future research should further investigate the benefits of
more integrated planning of patient transports within hospitals and outside of hospitals.

\subsection{Data availability}\label{sec:gap_data}
During our literature analysis, we found that very few data are currently available for future research. While most papers test their developed models with real data including real hospital layouts and distance matrices as well as real transport requests, others only use real layout data~\cite{vancroonenburg_2018a,nasir_2022a}. However, apart from~\cite{vancroonenburg_2018a}, who only use real layout data, none of the examined papers make their problem instances publicly available. The unavailability of real-world test instances hinders the development of new methods and the comparison of different modeling approaches and algorithms for patient transport problems. 

\subsection{Practical implementation}\label{sec:gap_practical}
Finally, we discuss possible reasons as to why, apart from~\cite{hanne_2009a}, none of the methods developed in the relevant literature have been implemented in daily hospital practice.

\medskip

Some of the papers state that no or only incomplete real data were available at the time of modeling~\cite{nickel_2000a,vancroonenburg_2018a,nasir_2022a}. One reason for this is that hospitals often do not systematically record data related to patient transports. As highlighted in~\cite{kergosien_2023a}, this results not only in difficulties in assessing the cost-benefit ratio of implementing OR/MS approaches, but also limits the ability to evaluate the quality of solutions obtained using these approaches. Therefore, ensuring a comprehensive data collection seems crucial in order to evaluate the potential practical impact of OR/MS approaches and foster their adoption in practice. Moreover, the collected data would facilitate the development of benchmark data sets and real-world test instances, which are currently missing from the patient transport literature as remarked in Section~\ref{sec:gap_data}.

\medskip

The most common reason for missing practical implementation mentioned, however, is a lack of uncertainty handling~\cite{melachrinoudis_2005a,beaudry_2010a,schmid_2014a,elmbach_2015a,elmbach_2019a}. Uncertainty handling, however, is very important since, as the papers also indicate, poor communication with regard to patient transport is normal, which often leads to uncertainty and incomplete information during patient transport planning. For instance, it is mentioned that upcoming transports are often not requested as soon as they become known, which would give some time for planning, but usually just before the desired pick-up time~\cite{beaudry_2010a}. This is also corroborated by the statement of our selected hospital not using any specialized software for patient transport planning, which has established the workflow of registering transports just in time via phone calls.

\medskip

Other reasons for missing or unsuccessful practical implementation mentioned in the literature are that hospital information systems are not considered suitable for the integration of the developed models~\cite{seguin_2019a}, not all desired constraints and special features could be modeled~\cite{elmbach_2019a,melachrinoudis_2005a}, the solvable instances are too small for practical use, or the optimization takes too long under real conditions~\cite{kallrath_2005a,turan_2011a,zhang_2015a}.

\medskip

All of these reasons have so far hindered practical implementation of patient transport models, which in turn results in a lack of experience with the cost-benefit ratio of implementation. In addition, many hospitals do not want to afford the high costs of commercial mixed-integer programming solvers~\cite{berg_2019a}. Overall, it appears that better modeling and handling of uncertainty, although already at least partially taken into account in various works, still represents a promising research area that is crucial for the practical implementation of the developed models. Moreover, taking all relevant practical constraints and the existing hospital information systems into account to integrate the developed models with these systems could help to foster the practical uptake of new approaches. 

\section{Conclusion}\label{sec:conclusion}
This paper reviews the existing OR/MS literature on intra-hospital patient transport and develops a five-field notation that allows a concise representation of different problem variants. This five-field notation can easily be extended (e.g., by adding further constraints or objectives) when new relevant aspects are introduced to the patient transport literature in the future.

\medskip

By analyzing the literature and comparing the results with insights from hospital practitioners, we identify important research gaps and directions for future work. These include improving the consideration of staff workload and ergonomic burden, addressing different qualification requirements of transport tasks, and exploring the integration of patient transport with material transport. Additionally, overcoming practical implementation challenges such as uncertainty handling, integration with hospital information systems, and long computation times is critical for enhancing the practical relevance of patient transport research.

\medskip

Moreover, the fact that patients need to be present for most tasks performed in a hospital underlines the importance of integrating patient transport with other operational planning problems in hospitals. This has so far only been considered explicitly in~\cite{schmid_2014a}, where patient transport is integrated with operating room scheduling. 

\medskip

Finally, although many of the papers examined test their models with fully or partially available real data, only~\cite{vancroonenburg_2018a} make their problem instances available for future research. Thus, creating an instance generator for different variants of patient transport problems could be useful to foster future research and to enable direct comparisons between different modeling and solution approaches.
\section*{Statements and declarations}

\subsection*{Declarations of interest}

\noindent
This research was partially funded by the Deutsche Forschungsgemeinschaft (DFG, German Research Foundation) -- Project number 443158418.

\subsection*{Author contributions}
\noindent
\textbf{Tom Lorenz Klein:} Validation, Formal Analysis, Investigation, Data Curation, Writing - Original Draft, Writing - Review and Editing, Visualization, Project Administration.
\textbf{Clemens Thielen:} Conceptualization, Methodology, Writing - Review and Editing, Supervision, Project Administration, Funding Acquisition.

\newpage

\bibliographystyle{elsarticle-num} 
\bibliography{cas-refs}

\newpage

\appendix

\section{Semi-structured interview guide}\label{sec:appa}

\subsection*{Introduction and background}
\begin{itemize}
    \item \textbf{Role and experience:} ``Can you please describe your role and experience in the patient transport department?''
    \item \textbf{Work experience:} ``How long have you been working in this area?''
    \item \textbf{Team size:} ``How many people work in the patient transport team?''
\end{itemize}

\subsection*{General organization of patient transport}
\begin{itemize}
    \item \textbf{Area of responsibility:} ``What tasks does patient transport involve in your hospital?''
    \item \textbf{Team structure:} ``Do you distinguish between medically trained and medically untrained transport staff in your team? If so, what is the practical impact of this distinction?''
    \item \textbf{Types of transport:} ``What types of transport requests does your team handle? Do you also transport materials?''
    \item \textbf{Organizational structure:} ``How is patient transport organized and structured in your hospital?''
\end{itemize}

\subsection*{Software and technology}
\begin{itemize}
    \item \textbf{Software use:} ``What software is used in the area of patient transport (if any)?''
    \item \textbf{Communication and requests:} ``How are transport requests communicated, and how are transporters informed about transport requests?''
    \item \textbf{Performance analysis:} ``Are evaluations carried out regarding the durations and punctuality of transports?''
\end{itemize}

\subsection*{Workload and employee management}
\begin{itemize}
    \item \textbf{Supporting aids:} ``Are pushing aids or other tools used to reduce the ergonomic burden of the transport staff?''
    \item \textbf{Workload management:} ``How is the workload distributed among employees, especially between male and female transport staff? Are there special regulations in this regard?''
\end{itemize}

\subsection*{Hospital to hospital transports and external service providers}
\begin{itemize}
    \item \textbf{Hospital to hospital transports:} ``How are patient transfers between different hospitals organized?''
    \item \textbf{Use of external service providers:} ``To what extent are external service providers used for transports?''
\end{itemize}

\subsection*{Discharges and planning}
\begin{itemize}
    \item \textbf{Patient discharges:} ``How are patient discharges organized and planned?''
    \item \textbf{Early notification of discharges:} ``Are patient discharges planned in advance and communicated early? How effective is this process?''
    \item \textbf{General planning of transports:} ``Are patient transports in your hospital generally planned and communicated well in advance, or do you frequently receive requests at very short notice? How does this affect the operational process?''
\end{itemize}

\subsection*{Challenges and possible solutions}
\begin{itemize}
    \item \textbf{Main challenges:} ``What are the biggest challenges in day-to-day patient transport operations?''
    \item \textbf{Coping strategies:} ``Have you implemented specific processes or techniques to overcome these challenges?''
    \item \textbf{Dealing with emergency situations:} ``How do you and your team deal with emergency situations?''
\end{itemize}

\subsection*{Conclusion}
\begin{itemize}
    \item \textbf{Final thoughts:} ``Are there any final thoughts or important aspects that we have not covered in this interview?''
\end{itemize}
\clearpage

\end{document}